\documentclass[12pt]{amsart}
\usepackage{amsmath, amssymb}
\usepackage{mathrsfs}
\newcommand{\no}[1]{#1}
\renewcommand{\no}[1]{}
\no{\usepackage{times}\usepackage[subscriptcorrection, slantedGreek, nofontinfo]{mtpro}
\renewcommand{\Delta}{\upDelta}}
\usepackage{color}
\usepackage{enumerate,comment}
\usepackage{amssymb}
\usepackage{graphicx}
\usepackage{dsfont}
\usepackage{epstopdf}
\graphicspath{{./pics/}}
\usepackage{booktabs}
\usepackage{threeparttable}

\date{\today}
\newcommand{\bel}{\begin{equation} \label}
\newcommand{\ee}{\end{equation}}
\def\beq{\begin{equation}}
\def\eeq{\end{equation}}
\newcommand{\bea}{\begin{eqnarray}}
\newcommand{\eea}{\end{eqnarray}}
\newcommand{\beas}{\begin{eqnarray*}}
\newcommand{\eeas}{\end{eqnarray*}}

\newcommand{\re}{\mathfrak R}

\newcommand{\R}{\mathbb{R}}
\newcommand{\C}{\mathbb{C}}

\newcommand{\pd}{\partial }
\newtheorem{theorem}{Theorem}[section]
\newtheorem{lemma}[theorem]{Lemma}

\newtheorem{defn}[theorem]{Definition}

\numberwithin{equation}{section}

\allowdisplaybreaks
\providecommand{\abs}[1]{\left\lvert#1\right\rvert}
\providecommand{\norm}[1]{\left\lVert#1\right\rVert}

\def\phi {\varphi}

\title[Inverse source problem with a posteriori boundary measurement]
{Inverse source problem with a posteriori boundary measurement for fractional diffusion equations}
\author[J. Janno]{Jaan Janno}
\address{Department of Cybernetics, Tallinn University of Technology
Ehitajate tee 5,19086 Tallinn, Estonia}
\email{ jaan.janno@ttu.ee}
\author[Y. Kian]{Yavar Kian}
\address{Aix Marseille Univ, Universit\'{e} de Toulon, CNRS, CPT, Marseille, France}
\email{yavar.kian@univ-amu.fr}

\date{}

\begin{document}

\begin{abstract}
In this article we study  inverse source problems for time-fractional diffusion equations from \textit{a posteriori} boundary measurement. Using the memory effect of these class of equations, we solve these  inverse problems for several class of space or time dependent source terms. We prove also the unique determination of a general class of space-time dependent separated variables source terms from such measurement. Our approach is based on the study of singularities of the Laplace transform in time of boundary traces of solutions of time-fractional diffusion equations.

\end{abstract}
\maketitle
\section{Introduction}
\subsection{Statement}
Let $\Omega\subset \mathbb{R}^d$ ($d\ge 2$) be an open bounded and connected subset with a
$ C^{2}$ boundary $\partial \Omega$. We define an elliptic operator $\mathcal{A}$  on the domain $\Omega$ by
\begin{equation}\label{A}
\mathcal{A} u(x) :=-\sum_{i,j=1}^d \partial_{x_i} \left( a_{i,j}(x) \partial_{x_j} u(x) \right)+q(x)u(x),\quad  x\in\Omega,
\end{equation}
where the potential $q \in L^\infty(\Omega)$ is nonnegative,
and the diffusion coefficient matrix $a:=(a_{i,j})_{1 \leq i,j \leq d} \in C^{1}
(\overline{\Omega};\mathbb{R}^{d\times d})$ is symmetric, i.e., $a_{i,j}(x)=a_{j,i}(x)$, for
any $ x \in \overline{\Omega}$, $i,j = 1,\ldots,d,$ and fulfills the following ellipticity condition
\begin{equation}\label{ell}
\exists c>0 :\,  \sum_{i,j=1}^d a_{i,j}(x) \xi_i \xi_j \geq c |\xi|^2,\quad  x \in \overline{\Omega},\ \xi=(\xi_1,\ldots,\xi_d) \in \mathbb{R}^d.
\end{equation}
Let the weight function $\rho\in L^\infty(\Omega)$ satisfy the condition
\begin{equation}\label{eqn:rho}
 0<c_0 \leq\rho(x) \leq C_0 <+\infty\quad\mbox{in } \Omega,
\end{equation}
with $c_0,C_0$ two positive constants. From now on and in all the remaining parts of this article we set $\R_+:=(0,+\infty)$ and $\mathbb N:=\{1,2,\ldots\}$. We introduce the Riemann-Liouville integral operator $I^\beta$ and the Riemann-Liouville fractional derivative  $D_t^\beta$, of order $\beta\in(0,1)$,  as follows
\[
I^\beta h(\cdot,t):=\frac{1}{\Gamma(\beta)}\int_0^t\frac{h(\cdot,\tau)}{(t-\tau)^{1-\beta}}\,d\tau,\quad D_t^\beta:=
\partial_t\circ I^{1-\beta}.
\]
We define also the Caputo fractional derivative  $\partial_t^\beta$, of order $\beta\in(0,1)$, by
$$\partial_t^\beta h=D_t^\beta (h-h(\cdot,0)),\quad h\in  C([0,+\infty);L^2(\Omega)).$$
We recall also that, for all $\beta\in(0,1)$, we have
$$\partial_t^\beta h:=I^{1-\beta}\partial_th,\quad h\in W^{1,1}_{loc}(\mathbb R_+;L^2(\Omega)).$$
Fixing $F\in L^1(\R_+;L^2(\Omega))$ and $\alpha \in (0, 1)$,
we consider weak solutions (in the sense of Definition \ref{d1}) $u$ of the following initial boundary value problem:
\begin{equation}\label{eq1}
\begin{cases}
\rho\partial_t^{\alpha}u +\mathcal{A} u =  F(x,t), & \mbox{in }\Omega\times\R_+,\\
 u= 0, & \mbox{on } \partial\Omega\times\R_+, \\
u=0, & \mbox{in } \Omega\times \{0\}.
\end{cases}
\end{equation}

With reference to \cite{K1,KY1}, one can check that \eqref{eq1} admits a unique weak solution lying in $L^1_{loc}(\R_+;H^s(\Omega))$, $1\leq s<2$.   Fixing $T\in\R_+$, in the present article,  we assume that the source term $F$ satisfies the following condition
\bel{source1} F(x,t)=\sigma(t)f(x),\quad t\in (0,T),\ x\in \Omega,\ee
with $\sigma\in L^1(0,T)$ and $f\in L^2(\Omega)$. We denote by $\pd_{\nu_a}$  the conormal derivative associated with the coefficient $a$ of $\pd\Omega$ defined by
$$\pd_{\nu_a}v(x)=\sum_{j=1}^da_{ij}(x)\partial_{x_j}v(x)\nu_i(x),\quad x\in\pd\Omega,$$
where $\nu=(\nu_1,\ldots,\nu_d)$ denotes the outward unit normal vector of $\pd\Omega$. We  assume that there exists $\delta\in(0,T)$ such that the time dependent part $\sigma$ of the source term in \eqref{eq1} satisfies the following condition
\bel{source} \sigma(t)=0,\quad t\in (T-\delta,T).\ee
Then, we consider three inverse problems:\\
(IP1)  Assuming that $\sigma$ is known and $\sigma\not\equiv0$, determine  $f$ from knowledge of  $\partial_{\nu_a} u(x,t)$, $(x,t)\in\Gamma\times(T-\epsilon,T)$ with $\Gamma$ an arbitrary open subset of $\partial\Omega$ and with $\epsilon\in(0,T)$ arbitrary small.\\
(IP2) Determine the space time dependent function $\sigma(t)f(x)$, $t\in(0,T)$ and $x\in\Omega$,  from knowledge of  $\partial_{\nu_a} u(x,t)$, $(x,t)\in\Gamma\times(T-\epsilon,T)$ with $\Gamma$ an arbitrary open subset of $\partial\Omega$ and with $\epsilon\in(0,T)$ arbitrary small.\\
(IP3) Assuming that $f$ is known and $f\not\equiv0$, determine $\sigma$ from the knowledge of $\partial_{\nu_a} u(x_0,t)$, $t\in(T-\epsilon,T)$ for some $x_0\in \pd\Omega$  and for $\epsilon\in(0,T)$ arbitrary small.\\

\subsection{Motivations}
We recall that systems of the form \eqref{eq1} model anomalous diffusion phenomena appearing  in applied sciences. This includes models of geophysics, environmental science and biology \cite{JR,NSY}. For such models, sub-diffusive  processes are described by \eqref{eq1}  and the kinetic equation \eqref{eq1} describes the corresponding macroscopic model to microscopic diffusion phenomena driven by continuous time random walk \cite{MK}.

In this context, the goal of the inverse problems (IP1)-(IP3) is to determine a source of anomalous diffusion from \textit{a posteriori} measurement. This class of problems may appear in several  practical situations such as environmental accident where the data are often  available only after the occurrence of such accident. We study this problem by exploiting the memory effect of solutions of \eqref{eq1}.

\subsection{Known results}

We recall that among the different formulations of inverse  problems, inverse source problems for time fractional diffusion equations  have received a lot of attention these last decades. The interested reader can refer to \cite{JR,LLY2} for an overview of this class of problems. Most of the results for this class of inverse problems are uniqueness results stated for source terms with separated variables of the form $F(x,t)=\sigma(t)f(x)$, $x\in\Omega$, $t\in(0,T)$. For this class of source terms, one can refer to  \cite{FK,LRY,LZ} for the determination of the time dependent component $\sigma(t)$ from measurement at one point and to \cite{JLLY,KSXY} for the determination of the space dependent component $f(x)$ from internal data. We refer also to \cite{KST} for the stable recovery of the space dependent component $f(x)$ in this class of problems. We can also mention the work of \cite{JKZ,KY2} dealing with the determination of a source term independent of one space direction in a cylindrical domain and the work of \cite{JK} who have considered the determination of a general space-time source term from the full knowledge of the solution close to the final time.

Most of the above mentioned results are stated with measurement throughout the full interval of time $(0,T)$. We are only aware of the three results \cite{JK,K2,KLY} considering measurement during an interval of time of the form $(T-\epsilon,T)$ with  $\epsilon\in(0,T)$ arbitrary small. In all these three works the authors use the memory effect of time-fractional diffusion equations exhibited by \cite[Theorem 1]{JK} which can not be applied in the context of the inverse problems (IP1)-(IP3).

\subsection{Main results}

Our main results will be devoted to the study of each of the inverse problems (IP1)-(IP3).

For (IP1) we obtain the following result.

\begin{theorem}\label{t1} Let $\sigma\in L^1(0,T)$, $f\in L^2(\Omega)$ with $\sigma\not\equiv0$ satisfying \eqref{source}. Consider also  $F\in L^1(\R_+;L^2(\Omega))$ satisfying \eqref{source1} and let $\Gamma$ be an open and not empty subset of $\partial\Omega$. Then for any $\epsilon\in(0,T)$ the following implication
\begin{equation}\label{t1b}\partial_{\nu_a}u_{|\Gamma\times (T-\epsilon,T)}\equiv0\Longrightarrow f\equiv0\end{equation}
holds true.
\end{theorem}

As a consequence of Theorem \ref{t1}, we obtain the following result for (IP2).

\begin{theorem}\label{t3} Let $\alpha$ be irrational and let  $\sigma_j\in L^1(0,T)$, $f_j\in L^2(\Omega)$, $j=1,2$. Assume also that condition \eqref{source} is fulfilled with $\sigma=\sigma_j$, $j=1,2$. For $j=1,2$, consider $F_j\in L^1(\R_+;L^2(\Omega))$, $j=1,2$, satisfying \eqref{source1} with $\sigma=\sigma_j$ and $f=f_j$, and let $u_j$ be the
weak solution of \eqref{eq1} corresponding to the source term $F=F_j$.
 Let also $\Gamma$ be an open and not empty subset of $\partial\Omega$. Then for any $\epsilon\in(0,T)$ the following implication
\begin{equation}\label{t3a}\partial_{\nu_a}u_{1}|_{\Gamma\times (T-\epsilon,T)}\equiv
\partial_{\nu_a}u_{2}|_{\Gamma\times (T-\epsilon,T)} \Longrightarrow F_{1}|_{\Omega\times (0,T)}\equiv F_{2}|_{\Omega\times (0,T)}\end{equation}
holds true.
\end{theorem}

Let us introduce the operator $A=\rho^{-1}\mathcal{A}$ acting on $L^2(\Omega;\rho d x)$ with domain $D(A):=\{h\in H^1_0(\Omega):\ \rho^{-1}\mathcal{A}\in L^2(\Omega;\rho d x)\}$. Here we denote by  $L^2(\Omega;\rho d x)$ the space
 of measurable functions $v$ satisfying
$\int_\Omega |v|^2\rho d x<\infty$ endowed with the inner product
$ \langle u,v\rangle =\int_\Omega uv\rho d x.$
Note that under condition \eqref{eqn:rho}, we have $L^2(\Omega;\rho d x)=L^2(\Omega)$
with equivalent norm, and thus we distinguish only the inner products but not the spaces. With reference to Section 2, we consider, for all $r>0$, $D(A^r)$ and we denote by $\rho D(A^r)$ the space of functions $h$ such that $\rho^{-1}h\in D(A^r)$. Assuming  that $\Omega$ has a $C^{2\lceil\frac{d}{4}\rceil+2}$  boundary, $a \in C^{1+2\lceil\frac{d}{4}\rceil}(\overline{\Omega};\R^{d^2})$, $\rho\in C^{1+2\lceil\frac{d}{4}\rceil}(\overline{\Omega})$ and applying the elliptic regularity of the operator $\rho^{-1}\mathcal{A}$ one can check that $D(A)=H^2(\Omega)\cap H^1_0(\Omega)$. In addition,  from \cite[Theorem 2.5.1.1]{Gr} we deduce  that, for all $\ell=1,\ldots,\lceil\frac{d}{4}\rceil+1$, we have
\begin{equation}\label{DH}D(A^\ell)=\{h\in H^{2\ell}(\Omega):\   h,\rho^{-1}\mathcal Ah,\ldots, (\rho^{-1}\mathcal A)^{\ell-1} h\in H^1_0(\Omega)\}.\end{equation}
By interpolation, for all $r>0$, $D(A^r)$ is embedded continuously into  $H^{2r}(\Omega)$. Combining this with the fact that $0$ is not in the spectrum of $A$ and applying the Sobolev embedding theorem, we deduce that for all $f\in \rho D(A^s)$,  with $s>\frac{d-2}{2}$, we have $A^{-k}\rho^{-1}f\in C^1(\overline{\Omega})$, $k\in\mathbb N$, and we can define the set
\begin{equation}\label{set}\mathcal G(f):= \{x\in\partial\Omega:\ \textrm{there exists $k\in\mathbb N\cup\{0\}$ such that } \partial_{\nu_a}A^{-k-2}\rho^{-1} f(x)\neq0\textrm{ and }\alpha(k+1)\not\in \mathbb N\}.\end{equation}
For (IP3) our  result is a uniqueness result which can be stated as follows.
\begin{theorem}
\label{t2} We assume that $\Omega$ has a $C^{2\lceil\frac{d}{4}\rceil+2}$  boundary, $a \in C^{1+2\lceil\frac{d}{4}\rceil}(\overline{\Omega};\R^{d^2})$, $\rho\in C^{1+2\lceil\frac{d}{4}\rceil}(\overline{\Omega})$.
Let  $\sigma\in L^1(0,T)$ satisfy the condition \eqref{source},  $f\in \rho D(A^s)$, with $s>\frac{d-2}{4}$, be a non-identically vanishing function and let $F\in L^1(\R_+;L^2(\Omega))$ satisfy \eqref{source1}. Then the  weak solution $u$ of \eqref{eq1}  is lying in $  L^1(0,T;C^1(\overline{\Omega}))$. Moreover, the set $\mathcal G(f)$, defined by \eqref{set}, is dense on $\partial\Omega$ and  for all  $x_0\in\mathcal G(f)$ and all $\epsilon\in(0,T)$ the implication
\begin{equation}\label{t2b}\left(\forall t\in(T-\epsilon,T),\ \partial_{\nu_a}u(x_0,t)=0\right)\Longrightarrow \sigma\equiv0\end{equation}
holds true. In addition, assuming that there exists $k_1\in\mathbb N\cup\{0\}$ and $g\in D(A^{k_1+s})$ of constant sign such that
$f=\rho A^{k_1}g$, we have $\mathcal G(f)=\partial\Omega$ and the implication \eqref{t2b} holds true for any $x_0\in\partial\Omega$.
\end{theorem}

Let us observe that the results of Theorem \ref{t1}, \ref{t3} and \ref{t2} are all stated with \textit{a posteriori} boundary measurement restricted to an arbitrary small interval of time of the form $(T-\epsilon,T)$ where $T$ denotes the final time. We are only aware of the two articles \cite{K2,KLY} studying this class of inverse source  problems with such  data. While the results of \cite{KLY} are stated with internal data, in the results of \cite{K2} the measurements are given by $\partial_{\nu_a} u(x,t)$ and $\partial_t^\alpha \partial_{\nu_a} u(x,t)$, $(x,t)\in\Gamma\times(T-\epsilon,T)$ with $\Gamma$ an arbitrary open subset of $\partial\Omega$ and with $\epsilon\in(0,T)$ arbitrary small. In both of these results, the authors apply the memory effect exhibited in \cite[Theorem 1]{JK} in order to determine from these class of data restricted to in interval of time of the form $(T-\epsilon,T)$ similar type of data on the full interval $(0,T)$. When the boundary measurement are only restricted to $\partial_{\nu_a} u(x,t)$, $(x,t)\in\Gamma\times(T-\epsilon,T)$, this approach does not work and a different strategy should be considered for the resolution of this problem. In the present article, we introduce a new strategy based on the study of the singularities of the Laplace transform in time of solutions of \eqref{eq1}, with  $F=0$ on $\Omega\times (T,+\infty)$, in order to analyze the inverse problems (IP1)-(IP3). Our results show that the memory effect of \eqref{eq1} can also be applied to the boundary measurement under consideration in inverse problems (IP1)-(IP3).

Note that in all the results of this article the source terms $F$ of \eqref{eq1} are unknown on $\Omega\times(T,+\infty)$ and we do not try to determine such values. This can be equivalently seen as  the statement of our inverse problems on the set $\Omega\times(0,T)$ instead of $\Omega\times\R_+$.

During the preparation of this article we realized that a result similar to Theorem \ref{t1} has been obtained, simultaneously and independently of this article, in \cite[Theorem 2]{Y} as a consequence of the main results of the same article. Our proof of Theorem \ref{t1} is completely different from the one under consideration in \cite{Y} where the author uses asymptotic properties of solutions while our analysis is mostly based on the study of singularities of Laplace transform in time of solution of \eqref{eq2} which coincides with the solution of \eqref{eq1} with  $F=0$ on $\Omega\times (T,+\infty)$.

Let us mention that, to the best of our knowledge, in Theorem \ref{t3} we obtain the first result of full determination of a general space-time dependent  source term of the form $\sigma(t)f(x)$ from boundary measurement. Indeed, in all the results that we are aware of (see e.g. \cite{JLLY,KLY,KSXY}) some \textit{a priori} knowledge of the time dependent component $\sigma(t)$ or  the space dependent component $f(x)$ of such source terms are required for the determination of this type of source terms. We recall also that there is a natural obstruction for the the determination of general source terms of the form $\sigma(t)f(x)$ for classical diffusion equations corresponding to \eqref{eq1} with $\alpha=1$ (see e.g. \cite[Theorem 7.1.]{KSXY}). In that sense, Theorem \ref{t3} emphasizes the particularity of time-fractional diffusion equations where in contrast to classical diffusion equations such inverse problems can be solved.

Notice that the result of Theorem \ref{t2} is stated with measurement at one point $x_0$ located on the explicit subset $\mathcal G(f)$ of $\partial\Omega$ defined by \eqref{set}. This formulation allows to state the result of Theorem \ref{t2} for general non-uniformly vanishing space dependent part $f(x)$ of
the source term $\sigma(t)f(x)$ while other similar results are often stated with some extra condition imposed to $f$ (see e.g. \cite{FK,LRY,LZ}). In addition, we prove in Theorem \ref{t2} that the set $\mathcal G(f)$ is dense in $\partial\Omega$ which means that the result of Theorem \ref{t2} holds true for at least one point in any arbitrary chosen open set of $\partial\Omega$. We exhibit also in Theorem \ref{t2} a general condition guarantying that $\mathcal G(f)=\partial\Omega$.

All the results of this article can be easily extended to \eqref{eq1} with homogeneous Neumann boundary condition instead of the homogeneous Dirichlet boundary condition. Moreover, with some minor modifications, our results can also be applied to super-diffusive models given by  \eqref{eq1} with $\alpha\in(1,2)$. Nevertheless, Theorem \ref{t3} is not true for classical diffusion equations  of the \eqref{eq1} with $\alpha=1$ and it is not clear whether
Theorem \ref{t1} and \ref{t2}  hold true  for classical diffusion equations.
\section{Preliminary properties}
In this section we introduce several preliminary properties related to the definition of solutions, the well-posedness and regularity properties for the initial boundary value problem \eqref{eq1}. We start by introducing the definition of weak solutions for problem \eqref{eq1} as follows.
\begin{defn}\label{d1}
We say that
$u\in L_{loc}^1(\mathbb R_+;L^2(\Omega))$ is a weak solution to \eqref{eq1} if it satisfies
the following conditions.
\begin{enumerate}
\item[{\rm(i)}] The identity $\rho(x)D^\alpha_t u(x,t) -\mathcal A u(x,t) = F(x,t),\  (x,t)\in   \Omega\times\R_+ $
 holds true in the sense of distributions in $\Omega\times\R_+$.
\item[{\rm(ii)}]  $I^{1-\alpha} u\in W_{loc}^{1,1}(\mathbb R_+;H^{-2}(\Omega))$ and $I^{1-\alpha} u(0, x)=0,\quad x\in\Omega$.
\item[{\rm(iii)}] $\tau_0=\inf\{\tau>0:\ e^{-\tau t}u\in L^1(\mathbb R_+;L^2(\Omega))\}<\infty$
and there exists $\tau_1\geq \tau_0$ such that for all $p\in\mathbb C$ satisfying $\mathfrak R\,p>\tau_1$ we have
$$
\hat{u}(\cdot,p):=\int_0^\infty e^{-p t}u(\cdot,t)\,d t\in H^1_0(\Omega).
$$
\end{enumerate}
\end{defn}

According to \cite{K1,KY1}, assuming that $F\in L^1(\R_+;L^2(\Omega))$, one can check that problem \eqref{eq1} admits a unique weak solution in the sense of Definition \ref{d1}.

Recall that the spectrum of the operator $A$
consists of an increasing sequence of strictly positive eigenvalues
$(\lambda_{n})_{n\geq1}$.  In the Hilbert space $L^2(\Omega;\rho dx)$, for each eigenvalue $\lambda_n$, we fix also $m_n\in\mathbb N$   the algebraic multiplicity of $\lambda_n$ and the family $\{\phi_{n,k}\}_{k=1}^{m_n}$ of eigenfunctions of $A$,
which forms an orthonormal basis in $L^2(\Omega;\rho dx)$ of the algebraic eigenspace of $A$ associated with $\lambda_n$. For all $s\geq 0$, we denote by $A^s$ the operator defined by
\[
A^s g=\sum_{n=1}^{+\infty}\sum_{k=1}^{m_n}\left\langle g,\phi_{n,k}\right\rangle
\lambda_{n}^s\phi_{n,k},\quad g\in D(A^s)
= \left\{h\in L^2(\Omega):\ \sum_{n=1}^{+\infty}\sum_{k=1}^{m_n}\abs{\left\langle g,
\phi_{n,k}\right\rangle}^2 \lambda_{n}^{2s}<\infty
\right\}
\]
and in $D(A^s)$ we introduce the norm
\[\|g\|_{D(A^s)}
= \left(\sum_{n=1}^{+\infty}\sum_{k=1}^{m_n}\abs{\left\langle g,
\phi_{n,k}\right\rangle}^2 \lambda_{n}^{2s}\right)^{\frac{1}{2}},
\quad g\in D(A^s).
\]
With reference to \cite{P}, we introduce the Mittag-Leffler function $E_{\alpha,\beta}(z)$ defined by

\begin{equation*}
  E_{\alpha,\beta}(z) = \sum_{n=0}^\infty \frac{z^n}{\Gamma(n\alpha+\beta)},\quad z\in \mathbb{C}
\end{equation*}
and we define the operator
\begin{equation}\label{S1}
S(t)h=\sum_{n=1}^\infty \sum_{k=1}^{m_n} t^{\alpha-1}E_{\alpha,\alpha}(-\lambda_n t^\alpha)\left\langle h,\phi_{n,k}\right\rangle\phi_{n,k},\quad h\in L^2(\Omega),\ t
\in \R_+.
\end{equation}
Let us consider the following initial boundary value problem
\begin{equation}\label{eq2}
\begin{cases}
\rho\partial_t^{\alpha}v +\mathcal{A} v =  \sigma(t)f(x), & \mbox{in }\Omega\times\R_+,\\
 v= 0, & \mbox{on } \partial\Omega\times\R_+, \\
v=0, & \mbox{in } \Omega\times \{0\},
\end{cases}
\end{equation}
where $\sigma$ denotes the extension of the function $\sigma\in L^1(0,T)$ by $0$ into an element of $L^1(\R_+)$.

It is well known  that the weak
solution of problem \eqref{eq2} is given by
$$v(\cdot,t)=\int_0^t\sigma(s)S(t-s)(\rho^{-1}f)ds,\quad t\in\R_+$$
and that $v=u$ on $\Omega\times (0,T)$, where $u$ denotes the solution of \eqref{eq1} with $F\in L^1(\R_+;L^2(\Omega))$ satisfying \eqref{source1}. In addition, following \cite[Proposition 2.2]{JK1}, \cite[Lemma 2.2]{JK2} and \cite[Proposition 2.1.]{KSXY}, one can check the following results

\begin{lemma}\label{l1} Let $f\in L^2(\Omega)$ and $\sigma\in L^1(0,T)$ satisfy \eqref{source}. Then  the problem \eqref{eq2} admits a unique weak solution $v\in L^1_{loc}(\R_+;H^{\frac{7}{4}}(\Omega))$ satisfying
\begin{equation}\label{l1a} \partial_{\nu_a}v(\cdot,t)|_{\Gamma}=\int_0^t\sigma(s) \partial_{\nu_a} [S(t-s)\rho^{-1}f]|_{\Gamma}ds,\quad t\in\R_+.\end{equation}
Moreover, the map $t\mapsto v(\cdot,t)$ is analytic with respect to $t\in(T-\delta,+\infty)$ as a map taking values in $H^{\frac{7}{4}}(\Omega)$. Finally, fixing $w=\partial_{\nu_a}v|_{\Gamma\times\R_+}$, we have
$$\inf\{p>0:\ e^{-pt}w(\cdot,t)\in L^1(\R_+;L^2(\Gamma))\}=0$$
and for all $p\in\mathbb C_+:=\{z\in\mathbb C:\ \re(z)>0\}$ we have
\begin{equation}\label{l1b}\partial_{\nu_a}\hat{v}(\cdot,p)|_{\Gamma}=\hat{w}(\cdot,p)=\hat{\sigma}(p)\sum_{n=1}^\infty\sum_{k=1}^{m_n} \frac{\left\langle \rho^{-1}f,\phi_{n,k}\right\rangle}{\lambda_n+p^\alpha}\partial_{\nu_a}\phi_{n,k}|_{\Gamma}.\end{equation}

\end{lemma}

\begin{lemma}\label{l2} Assume that $\Omega$ has a $C^{2\lceil\frac{d}{4}\rceil+2}$  boundary, $a \in C^{1+2\lceil\frac{d}{4}\rceil}(\overline{\Omega};\R^{d^2})$, $\rho\in C^{1+2\lceil\frac{d}{4}\rceil}(\overline{\Omega})$.  Let $f\in \rho D(A^s)$, with $s>\frac{d-2}{4}$ and let $\sigma\in L^1(0,T)$  satisfy \eqref{source}. Then the problem \eqref{eq2} admits a unique weak solution $v\in L^1_{loc}(\R_+;C^1(\overline{\Omega}))$ satisfying
\begin{equation}\label{l2a} \partial_{\nu_a}v(x,t)=\int_0^t\sigma(s) \partial_{\nu_a} [S(t-s)\rho^{-1}f](x)ds,\quad t\in(0,+\infty),\ x\in\partial\Omega.\end{equation}
 Moreover, the map $t\mapsto v(\cdot,t)$ is analytic with respect to $t\in(T-\delta,+\infty)$ as a map taking values  in $C^1(\overline{\Omega})$. Finally, for any $x_1\in\partial\Omega$  fixing $w_1:=\R_+\ni t\mapsto \partial_{\nu_a} v(x_1,t)$ we have
$$\inf\{p>0:\ e^{-pt}w_1(t)\in L^1(\R_+)\}=0$$
and for all $p\in\mathbb C_+$ we have
\begin{equation}\label{l2b}\hat{w_1}(p):=\hat{\sigma}(p)\sum_{n=1}^\infty\sum_{k=1}^{m_n} \frac{\left\langle \rho^{-1}f,\phi_{n,k}\right\rangle\partial_{\nu_a}\phi_{n,k}(x_1)}{\lambda_n+p^\alpha}.\end{equation}
\end{lemma}

\section{Proof of Theorem \ref{t1}}

 Let $\epsilon\in(0,T)$ and let $\Gamma$ be an open and not empty subset of $\partial\Omega$. Assume that the solution $u$ of \eqref{eq1} satisfies $\partial_{\nu_a}u_{|\Gamma\times (T-\epsilon,T)}\equiv0$. Then, recalling that the solution $v$ of \eqref{eq2} satisfies $v=u$ on $\Omega\times (0,T)$, we deduce that $\partial_{\nu_a}v_{|\Gamma\times (T-\epsilon,T)}=\partial_{\nu_a}u_{|\Gamma\times (T-\epsilon,T)}\equiv0$. Without loss of generality, we assume that $\epsilon=\delta$ with $\delta>0$ appearing in \eqref{source}.  Applying Lemma \ref{l1}, we deduce that $(T-\delta,+\infty)\ni t\mapsto\partial_{\nu_a}v(\cdot,t)_{|\Gamma}$ is an analytic function taking values in $L^2(\Gamma)$. Therefore, fixing $w=\partial_{\nu_a}v_{|\Gamma\times \R_+}$ and applying unique continuation for holomorphic functions, we find
\bel{t1a}w(\cdot,t)=\partial_{\nu_a}v(\cdot,t)_{|\Gamma}\equiv 0,\quad t\in(T-\epsilon,+\infty).\ee
Thus the Laplace transform in time $\hat{w}(\cdot,p)$ of $w$ is holomorphic with respect  to $p\in\mathbb C$ as a function taking values in $L^2(\Gamma)$. On the other hand, applying \eqref{l1b}, we deduce that the following identity
\bel{t1b}\hat{w}(\cdot,p)=\hat{\sigma}(p)\sum_{n=1}^\infty \frac{\sum_{k=1}^{m_n}\left\langle \rho^{-1}f,\phi_{n,k}\right\rangle\partial_{\nu_a}\phi_{n,k}|_{\Gamma}}{\lambda_n+p^\alpha}\ee
holds true for all $p\in\mathbb C_+$.
In the same way, one can easily check that the map
$$\mathcal K:=p\mapsto \sum_{n=1}^\infty \frac{\sum_{k=1}^{m_n}\left\langle \rho^{-1}f,\phi_{n,k}\right\rangle}{\lambda_n+p^\alpha}\partial_{\nu_a}\phi_{n,k}|_{\Gamma}$$
admits an holomorphic extension to $\mathbb C\setminus(-\infty,0] $ as a map taking values in $L^2(\Gamma)$. Here and in all the remaining parts of the article, for any $p\in\mathbb C\setminus(-\infty,0]$ we set $p^\alpha=e^{\alpha \log(p)}$ with $\log$ the complex logarithm defined on  $\mathbb C\setminus(-\infty,0]$. Therefore, the identity \eqref{t1b} holds true for all $p\in\mathbb C\setminus(-\infty,0]$. Now let us fix $R>0$ and consider the identity \eqref{t1b} for $p=Re^{i\theta}$ with $\theta\in(-\pi,\pi)$. Sending $\theta\to\pm\pi$ and using the fact that $\hat{w}(\cdot,p)$ and $\hat{\sigma}(p)$ are holomorphic with respect  to $p\in\mathbb C$, we obtain
$$\hat{w}(\cdot,-R)=\hat{w}(\cdot,Re^{\pm i\pi})=\hat{\sigma}(Re^{\pm i\pi})\sum_{n=1}^\infty \frac{\sum_{k=1}^{m_n}\left\langle \rho^{-1}f,\phi_{n,k}\right\rangle\partial_{\nu_a}\phi_{n,k}|_{\Gamma}}{\lambda_n+R^\alpha e^{\pm i\alpha\pi}}.$$
Taking the difference of these two expressions, we obtain
$$\begin{aligned}&\hat{\sigma}(-R)\sum_{n=1}^\infty \left(\sum_{k=1}^{m_n}\left\langle \rho^{-1}f,\phi_{n,k}\right\rangle\partial_{\nu_a}\phi_{n,k}|_{\Gamma}\right)\left(\frac{1}{\lambda_n+R^\alpha e^{ i\alpha\pi}}-\frac{1}{\lambda_n+R^\alpha e^{ -i\alpha\pi}}\right)\\
&=\hat{w}(\cdot,-R)-\hat{w}(\cdot,-R)\equiv0.\end{aligned}$$
It follows that
$$\begin{aligned}0&\equiv\hat{\sigma}(-R)\sum_{n=1}^\infty \left(\sum_{k=1}^{m_n}\left\langle \rho^{-1}f,\phi_{n,k}\right\rangle\partial_{\nu_a}\phi_{n,k}|_{\Gamma}\right)\left(\frac{R^\alpha (e^{ -i\alpha\pi}- e^{ i\alpha\pi})}{(\lambda_n+R^\alpha e^{ i\alpha\pi})(\lambda_n+R^\alpha e^{ -i\alpha\pi})}\right)\\
&=\hat{\sigma}(-R)\sum_{n=1}^\infty \left(\sum_{k=1}^{m_n}\left\langle \rho^{-1}f,\phi_{n,k}\right\rangle\partial_{\nu_a}\phi_{n,k}|_{\Gamma}\right)\left(\frac{-2iR^\alpha \sin(\alpha\pi)}{(\lambda_n+R^\alpha e^{ i\alpha\pi})(\lambda_n+e^{ -2i\alpha\pi}R^\alpha e^{ i\alpha\pi} )}\right).\end{aligned}$$
Recalling that $\alpha\in(0,1)$, we find $\sin(\alpha\pi)\neq0$ and we deduce that
$$\hat{\sigma}(-R)\sum_{n=1}^\infty \left(\frac{\sum_{k=1}^{m_n}\left\langle \rho^{-1}f,\phi_{n,k}\right\rangle\partial_{\nu_a}\phi_{n,k}|_{\Gamma}}{(\lambda_n+R^\alpha e^{ i\alpha\pi})(\lambda_n+e^{ -2i\alpha\pi}R^\alpha e^{ i\alpha\pi} )}\right)\equiv 0.$$
In the same way, using the fact that $\hat{\sigma}$ is holomorphic and $\sigma\not\equiv0$, we can find $R_1,R_2\in\R_+$ with $R_1<R_2$ such that $|\hat{\sigma}(-R)|>0$ for $R\in(R_1,R_2)$. Then, it follows that
\bel{t1c}\sum_{n=1}^\infty \left(\frac{\sum_{k=1}^{m_n}\left\langle \rho^{-1}f,\phi_{n,k}\right\rangle\partial_{\nu_a}\phi_{n,k}|_{\Gamma}}{(\lambda_n+r e^{ i\alpha\pi})(\lambda_n+e^{ -2i\alpha\pi}r e^{ i\alpha\pi} )}\right)\equiv 0,\quad r\in(R_1^\alpha,R_2^\alpha).\ee
Now, let us consider the set $\mathcal O:=\mathbb C\setminus \{-\lambda_n,-e^{ -2i\alpha\pi}\lambda_n:\ n\in\mathbb N\}$ and the
map
$$\mathcal H:\mathcal O\ni z\mapsto \sum_{n=1}^\infty \left(\frac{\sum_{k=1}^{m_n}\left\langle \rho^{-1}f,\phi_{n,k}\right\rangle\partial_{\nu_a}\phi_{n,k}|_{\Gamma}}{(\lambda_n+z)(\lambda_n+e^{ -2i\alpha\pi}z )}\right).$$
Since $f\in L^2(\Omega)$, one can easily check that the map $\mathcal H$ is an holomorphic function on $\mathcal O$ as a map taking values in $L^2(\Gamma)$. Combining this with \eqref{t1c} and applying the unique continuation for holomorphic functions, we obtain the following identity
\bel{t1d}\sum_{n=1}^\infty \left(\frac{\sum_{k=1}^{m_n}\left\langle \rho^{-1}f,\phi_{n,k}\right\rangle\partial_{\nu_a}\phi_{n,k}|_{\Gamma}}{(\lambda_n+z)(\lambda_n+e^{ -2i\alpha\pi}z )}\right)\equiv 0,\quad z\in\mathcal O.\ee
Fixing $n\in\mathbb N$ and multiplying \eqref{t1d} by $(z+\lambda_n)$ and sending $z\to-\lambda_n$, we find
\bel{t1e}\sum_{k=1}^{m_n}\left\langle \rho^{-1}f,\phi_{n,k}\right\rangle\partial_{\nu_a}\phi_{n,k}|_{\Gamma}\equiv0,\quad n\in\mathbb N.\ee
On the other hand, it is well known that the maps $\partial_{\nu_a}\phi_{n,k}|_{\Gamma}$, $k=1,\ldots,m_n$ are linearly independent (see e.g. \cite[Lemma 2.1.]{CK} or \cite[Step 4 in the proof of Theorem 1.1]{KSXY} for the proof of an equivalent property) and applying \eqref{t1e}, we find
$$\left\langle \rho^{-1}f,\phi_{n,k}\right\rangle=0,\quad n\in\mathbb N,\ k=1,\ldots,m_n.$$
From this last identity, we deduce that $\rho^{-1}f\equiv0$ and \eqref{eqn:rho} implies that $f\equiv0$. This completes the proof of the theorem.
\section{Proof of Theorem \ref{t3}}

Let us assume that  \bel{t3b}\partial_{\nu_a}u_{1|\Gamma\times (T-\epsilon,T)}\equiv
\partial_{\nu_a}u_{2|\Gamma\times (T-\epsilon,T)}.\ee  Let $v_j$, $j=1,2,$ be the solutions of
\eqref{eq2} that correspond to the pairs $(f,\sigma)=(f_j,\sigma_j)$, $j=1,2$, respectively and notice that \eqref{t3b} and Lemma \ref{l1} imply that
\bel{t3c}\partial_{\nu_a}v_{1}|_{\Gamma\times (T-\epsilon,+\infty)}\equiv
\partial_{\nu_a}v_{2}|_{\Gamma\times (T-\epsilon,+\infty)}.\ee

Further, let us define $w=\partial_{\nu_a}v_1(\cdot,t)_{|\Gamma} - \partial_{\nu_a}v_2(\cdot,t)_{|\Gamma}$, $t>0$. In view of \eqref{t3c} and condition \eqref{source} with $\sigma=\sigma_j$, $j=1,2$, the Laplace transforms of $w$, $\sigma_1$ and $\sigma_2$ are holomorphic in $\C$ and
 the formula
\bel{theorF2}\hat{w}(\cdot,p)=\hat{\sigma_1}(p)\sum_{n=1}^\infty
\frac{\sum_{k=1}^{m_n}\left\langle \rho^{-1}f_1,\phi_{n,k}\right\rangle\partial_{\nu_a}\phi_{n,k}|_{\Gamma}}{\lambda_n+p^\alpha}
-\hat{\sigma_2}(p)\sum_{n=1}^\infty
\frac{\sum_{k=1}^{m_n}\left\langle \rho^{-1}f_2,\phi_{n,k}\right\rangle\partial_{\nu_a}\phi_{n,k}|_{\Gamma}}{\lambda_n+p^\alpha}\ee
is valid for $p\in\C\setminus (-\infty,0]$.
We fix $R>0$ and we set $p=Re^{i\theta}$ with $\theta\in(-\pi,\pi)$ in \eqref{theorF2}. Sending
$\theta\to\pm\pi$  and taking the difference of expressions obtained, we get
$$\hat{\sigma}_1(-R)\sum_{n=1}^\infty \left(\frac{\sum_{k=1}^{m_n}\left\langle
\rho^{-1}f_1,\phi_{n,k}\right\rangle\partial_{\nu_a}\phi_{n,k}|_{\Gamma}}
{(\lambda_n+R^\alpha e^{ i\alpha\pi})(\lambda_n+R^\alpha e^{- i\alpha\pi} )}\right)-
\hat{\sigma}_2(-R)\sum_{n=1}^\infty \left(\frac{\sum_{k=1}^{m_n}\left\langle
\rho^{-1}f_2,\phi_{n,k}\right\rangle\partial_{\nu_a}\phi_{n,k}|_{\Gamma}}
{(\lambda_n+R^\alpha e^{ i\alpha\pi})(\lambda_n+R^\alpha e^{- i\alpha\pi} )}\right)\equiv 0,$$
where $R>0$. Replacing $R^\alpha$ by $r>0$, for all $r\in\R_+$, we have
$$\hat{\sigma}_1(-r^{1\over\alpha})\sum_{n=1}^\infty \left(\frac{\sum_{k=1}^{m_n}\left\langle
\rho^{-1}f_1,\phi_{n,k}\right\rangle\partial_{\nu_a}\phi_{n,k}|_{\Gamma}}
{(\lambda_n+r e^{ i\alpha\pi})(\lambda_n+r e^{ -i\alpha\pi} )}\right)-
\hat{\sigma}_2(-r^{1\over\alpha})\sum_{n=1}^\infty \left(\frac{\sum_{k=1}^{m_n}\left\langle
\rho^{-1}f_2,\phi_{n,k}\right\rangle\partial_{\nu_a}\phi_{n,k}|_{\Gamma}}
{(\lambda_n+r e^{ i\alpha\pi})(\lambda_n+r e^{- i\alpha\pi} )}\right)\equiv 0.$$

Let us define the set $\mathcal O_1:=\mathbb C\setminus( \{-e^{ \pm i\alpha\pi}\lambda_n:\ n\in\mathbb N\}\cup (-\infty,0])$ and
the
map
$$\mathcal H_1:\mathcal O_1\ni
 z\mapsto \hat{\sigma}_1(-z^{1\over\alpha})\Psi_1(\cdot,z)+
\hat{\sigma}_2(-z^{1\over\alpha})\Psi_2(\cdot,z),$$
with
$$
\Psi_1(\cdot,z)=\sum_{n=1}^\infty \left(\frac{\sum_{k=1}^{m_n}\left\langle \rho^{-1}f_1,
 \phi_{n,k}\right\rangle\partial_{\nu_a}\phi_{n,k}|_{\Gamma}}{(\lambda_n+z e^{ i\alpha\pi})(\lambda_n+z e^{- i\alpha\pi} )}\right),$$
$$\Psi_2(\cdot,z)=-\sum_{n=1}^\infty \left(\frac{\sum_{k=1}^{m_n}\left\langle \rho^{-1}f_2,
 \phi_{n,k}\right\rangle\partial_{\nu_a}\phi_{n,k}|_{\Gamma}}{(\lambda_n+z e^{ i\alpha\pi})(\lambda_n+z e^{- i\alpha\pi} )}\right).
$$
Then, for all $r>0$, we have $\mathcal H_1(r)=0$. Moreover, in a similar way to Theorem \ref{t1}, one can check
 that the map $\mathcal H_1$ is holomorphic in $\mathcal O_1$. Therefore, by unique continuation of holomorphic functions
  $\mathcal H_1$ vanishes  in $\mathcal O_1$. Now
we set $z=\varrho e^{i\theta}$ in the relation $\mathcal H_1(z)\equiv 0$  and send $\theta$ to $\pm\pi$.
We obtain the following system of equations:
\bea\label{theorF3}
&&\hat{\sigma}_1(-\varrho^{1\over\alpha}e^{i\pi\over\alpha})\Psi_1(\cdot,-\varrho)+
\hat{\sigma}_2(-\varrho^{1\over\alpha}e^{i\pi\over\alpha})\Psi_2(\cdot,-\varrho)\equiv 0,
\\ \label{theorF4}
&&\hat{\sigma}_1(-\varrho^{1\over\alpha}e^{-{i\pi\over\alpha}})\Psi_1(\cdot,-\varrho)+
\hat{\sigma}_2(-\varrho^{1\over\alpha}e^{-i{\pi\over\alpha}})\Psi_2(\cdot,-\varrho)\equiv 0,
\eea
where $\varrho>0$. Since the maps $\Psi_j$ are holomorphic in
  the set $\mathcal O_2=\C\setminus  \{-e^{ \pm i\alpha\pi}\lambda_n:\ n\in\mathbb N\}$, either
\\
(1) $\Psi_1=\Psi_2\equiv 0,\;\; j=1,2$
\\
or\\
(2) there exist $ j_0\in \{1;2\}$ and  $r_0>0$ such that $\norm{\Psi_{j_0}(\cdot,z)}_{L^2(\Gamma)}> 0$ for $0<|z|<r_0$.
\\
In case (1), using the method presented in the proof of Theorem \ref{t1} we deduce the relations $f_1\equiv 0$ and $f_2\equiv 0$ that imply
 $F_1|_{\Omega\times(0,T)}= F_2|_{\Omega\times(0,T)}\equiv 0$. In case (2) the determinant of the system \eqref{theorF3}, \eqref{theorF4}
  is zero for $0<\varrho<r_0$, i.e.,
\bea
\label{theorF5} \det \left(\begin{array}{lll}
&\hat{\sigma}_1(-\varrho^{1\over\alpha}e^{i\pi\over\alpha}) &\hat{\sigma}_2(-\varrho^{1\over\alpha}e^{i\pi\over\alpha})
\\
&\hat{\sigma}_1(-\varrho^{1\over\alpha}e^{-{i\pi\over\alpha}}) &\hat{\sigma}_2(-\varrho^{1\over\alpha}e^{-{i\pi\over\alpha}})
\end{array}\right)= 0,\;\; 0<\varrho<r_0.
\eea
If $\hat{\sigma}_1\equiv0$, combining \eqref{theorF3} with the arguments used at the end of the proof of Theorem \ref{t1}
 we deduce that either $\sigma_2\equiv0$ or $f_2\equiv0$ and it follows that $F_1|_{\Omega\times(0,T)}= F_2|_{\Omega\times(0,T)}\equiv 0$.
  In the same way, assuming that $\hat{\sigma}_2\equiv0$ we deduce that $F_1|_{\Omega\times(0,T)}= F_2|_{\Omega\times(0,T)}\equiv 0$.
  Therefore, we need to consider the situation where (2) holds true and $\hat{\sigma}_j\not\equiv0$, $j=1,2$.
  Since $\hat{\sigma}_1$ and $\hat{\sigma}_2$ are entire and not uniformly vanishing, we can find $r_1\in(0,r_0)$ such that
$$\min(|\hat{\sigma}_1(z)|,|\hat{\sigma}_2(z)|)>0,\quad 0<|z|<r_1.$$
In that situation, either $\hat{\sigma}_1(z)/\hat{\sigma}_2(z)$ or
$\hat{\sigma}_2(z)/\hat{\sigma}_1(z)$ will have a finite limit as $z\to0$ and then can be extended to an holomorphic
 function with respect to  $z\in D_{r_1}:=\{\eta\in\mathbb C:\ |\eta|<r_1\}$. From now on, without loss of generality,
  we assume that $\hat{\sigma}_2(z)/\hat{\sigma}_1(z)$ admits an holomorphic extension with respect to  $z\in D_{r_1}$.
   Then there exists a function $\mathcal C$ holomorphic
in $D_{r_1}$  such that $$\hat{\sigma}_2(z)= \mathcal C(z)\hat{\sigma}_1(z),\quad z\in D_{r_1}.$$
 Using this notation, we can transform the equation \eqref{theorF5} into
$$
\hat{\sigma}_1(-\varrho^{1\over\alpha}e^{i\pi\over\alpha})\hat{\sigma}_1(-\varrho^{1\over\alpha}e^{-{i\pi\over\alpha}})
\left(\mathcal C(-\varrho^{1\over\alpha}e^{i\pi\over\alpha})-\mathcal C(-\varrho^{1\over\alpha}e^{-{i\pi\over\alpha}})\right)= 0,
\;\;0<\varrho<r_1.
$$
This implies that
 $\mathcal C(-\varrho^{1\over\alpha}e^{i\pi\over\alpha})\equiv \mathcal C(-\varrho^{1\over\alpha}e^{-{i\pi\over\alpha}})$ for $0<\varrho<r_1$.
\\
 Due to the holomorphy of $\mathcal C$, we deduce that
\bea\label{theorF6}
\mathcal C(z)=\mathcal C(ze^{-{2i\pi\over\alpha}}),\;\; |z|<r_1^{1\over\alpha}.
\eea
Differentiating this relation at $z=0$ we have
$$
\mathcal C^{(n)}(0)=\mathcal C^{(n)}(0)e^{-{2ni\pi\over\alpha}},\;\; n\in\mathbb N.
$$
Since $\alpha$ is irrational, we have $e^{-{2ni\pi\over\alpha}}\ne 1$, $n\in\mathbb N$. Consequently, $\mathcal C^{(n)}(0)=0$, $n\in\mathbb N$.
Hence, the function $\mathcal C$ is constant and we have $\sigma_2\equiv \mathcal C\sigma_1$. Using the fact that $\hat{\sigma}_j\not\equiv0$, $j=1,2$,
we deduce that $\mathcal C\neq 0$. Then, we can rewrite $\sigma_2(t)f_2(x)$ in the form $$\sigma_2(t)f_2(x)=\sigma_1(t)\tilde f_2(x),\quad (x,t)\in\Omega\times(0,T),$$ where $\tilde f_2\equiv \mathcal C f_2$.
Then, the function $u=v_1-v_2$ solves \eqref{eq1} with the source function
$F(x,t)=(f_1(x)-\tilde f_2(x))\sigma_1(t)$ and satisfies the condition
$\partial_{\nu_a}u_{|\Gamma\times (T-\epsilon,T)}\equiv 0$. Thus,
Theorem \ref{t1} implies $f_1-\tilde f_2\equiv 0$. Consequently, we have
$$\sigma_2(t)f_2(x)=\sigma_1(t)\mathcal C f_2(x)=\sigma_1(t)f_1(x),\quad (x,t)\in\Omega\times(0,T)$$
which implies  that $F_1|_{\Omega\times(0,T)}\equiv F_2|_{\Omega\times(0,T)}$.
 This completes the proof of Theorem \ref{t3}.

\section{Proof of Theorem \ref{t2}}
In this section we use the notation of Section 2. Recall that  $\rho^{-1}f\in D(A^{s})$, with $s>\frac{d-2}{4}$, and we deduce that $(A+re^{i\alpha\pi})^{-1}(A+re^{-i\alpha\pi})^{-1}\rho^{-1}f\in H^{2s+2}(\Omega)$, $r>0$, and the Sobolev embedding theorem implies that $(A+re^{i\alpha\pi})^{-1}(A+re^{-i\alpha\pi})^{-1}f\in C^1(\overline{\Omega})$. Therefore, we can define the set $\mathcal J(f)$ of points $x_0\in\partial\Omega$ such that for all  $\epsilon>0$ there exists  $r_0\in(0,\epsilon)$ such that the condition
\begin{equation}\label{cond2}\partial_{\nu_a}(A+r_0e^{i\alpha\pi})^{-1}(A+r_0e^{-i\alpha\pi})^{-1}\rho^{-1}f(x_0)\neq0\end{equation}
is fulfilled. The proof of Theorem \ref{t2} will be divided into four steps. We will start by proving that for $\mathcal G(f)$,  defined by \eqref{set}, we have  $\mathcal G(f)=\mathcal J(f)$. Then, using this identity, we will prove that the set $\mathcal G(f)$ is dense in $\partial\Omega$. Then, we prove that, for all $x_0\in\mathcal G(f)$ the implication \eqref{t2b} holds true.  Finally, we will show that when $f=\rho A^{k_1}g$, with $g\in D(A^{k_1+s})$ of constant sign, we have $\mathcal G(f)=\partial\Omega$.

\textbf{Step 1.} In this step we will show that $\mathcal J(f)=\mathcal G(f)$. For this purpose, it would be enough to prove that $\partial\Omega\setminus\mathcal G(f)=\partial\Omega\setminus\mathcal J(f)$. Let us fix
$\delta_1=\left(2\norm{A^{-1}}_{\mathcal B(D(A^{s+1}))}\right)^{-1}$. Since $\rho^{-1}f\in D(A^s)$, we deduce that  the map $r\mapsto (A+re^{i\alpha\pi})^{-1}(A+re^{-i\alpha\pi})^{-1}\rho^{-1}f$ is analytic with respect to $r\in(0,+\infty)$ as a function taking values in $D(A^{s+1})$. Combining this with the continuous embedding $D(A^{s+1})\subset C^1(\overline{\Omega})$, one can check  that the map $r\mapsto (A+re^{i\alpha\pi})^{-1}(A+re^{-i\alpha\pi})^{-1}\rho^{-1}f$ is analytic with respect to $r\in(0,+\infty)$ as a function taking values in $C^1(\overline{\Omega})$. From this last property and the unique continuation for analytic functions, we deduce that   $x\in \partial\Omega\setminus\mathcal J(f)$ if and only if $x\in\partial\Omega$ satisfies the following condition
\begin{equation}\label{cond3}\exists\epsilon_1\in(0,\delta_1),\quad  \partial_{\nu_a}(A+re^{i\alpha\pi})^{-1}(A+re^{-i\alpha\pi})^{-1}\rho^{-1}f(x)=0,\quad r\in(0,\epsilon_1) .\end{equation}
On the other hand, for $z_1,z_2\in D_{\delta_1}:=\{\eta\in\mathbb C:\ |\eta|<\delta_1\}$, we have
$$\begin{aligned}(A+z_1e^{i\alpha\pi})^{-1}(A+z_2e^{-i\alpha\pi})^{-1}\rho^{-1}f&=(Id+z_1e^{i\alpha\pi}A^{-1})^{-1}(Id+z_2e^{-i\alpha\pi}A^{-1})^{-1}A^{-2}\rho^{-1}f\\
&=\sum_{j=0}^\infty\sum_{k=0}^\infty (-e^{i\alpha\pi})^{k}z_1^k(-e^{-i\alpha\pi})^{j}z_2^jA^{-j-k-2}\rho^{-1}f.\end{aligned}$$
Using the fact $\rho^{-1}f\in D(A^s)$ and applying the Sobolev embedding theorem, we deduce that, for all $z_1,z_2\in D_{\delta_1}$, we have
$$\begin{aligned}&\sum_{j=0}^\infty\sum_{k=0}^\infty \norm{(-e^{i\alpha\pi})^{k}z_1^k(-e^{-i\alpha\pi})^{j}z_2^jA^{-j-k-2}\rho^{-1}f}_{C^1(\overline{\Omega})}\\
&\leq C\sum_{j=0}^\infty\sum_{k=0}^\infty \norm{(-e^{i\alpha\pi})^{k}z_1^k(-e^{-i\alpha\pi})^{j}z_2^jA^{-j-k-2}\rho^{-1}f}_{H^{2s+2}(\Omega)}\\
\ &\leq  C\sum_{j=0}^\infty\sum_{k=0}^\infty \norm{(-e^{i\alpha\pi})^{k}z_1^k(-e^{-i\alpha\pi})^{j}z_2^jA^{-j-k-2}\rho^{-1}f}_{D(A^{s+1})}\\
\ &\leq C\sum_{j=0}^\infty\sum_{k=0}^\infty \delta_1^j\norm{A^{-1}}_{\mathcal B(D(A^{s+1}))}^{j+1}\delta_1^k\norm{A^{-1}}_{\mathcal B(D(A^{s+1}))}^{k}\norm{A^{-1}}_{\mathcal B(D(A^{s+1}),D(A^{s}))}\norm{\rho^{-1}f}_{ D(A^s)}\\
\ &\leq C\norm{A^{-1}}_{\mathcal B(D(A^{s+1}),D(A^{s}))}\norm{\rho^{-1}f}_{ D(A^s)}\sum_{j=0}^\infty\sum_{k=0}^\infty 2^{-k}2^{-j}<\infty.\end{aligned}$$

Therefore, the sequence
$$\sum_{j=0}^{N_1}\sum_{k=0}^{N_2} (-e^{i\alpha\pi})^{k}z_1^k(-e^{-i\alpha\pi})^{j}z_2^jA^{-j-k-2}\rho^{-1}f,\quad N_1,\ N_2\in\mathbb N$$
converges uniformly with respect to $z_1,z_2\in D_{\delta_1}$ in the sense of  $C^1(\overline{\Omega})$ valued function and, for any $x\in\partial\Omega$,  we have
\begin{equation}\label{cond4}\begin{aligned}&\partial_{\nu_a}(A+z_1e^{i\alpha\pi})^{-1}(A+z_2e^{-i\alpha\pi})^{-1}\rho^{-1}f(x)\\
&=\sum_{j=0}^\infty\sum_{k=0}^\infty (-e^{i\alpha\pi})^{k}z_1^k(-e^{-i\alpha\pi})^{j}z_2^j\partial_{\nu_a}A^{-j-k-2}\rho^{-1}f(x),\quad z_1,z_2\in D_{\delta_1}.\end{aligned}\end{equation}
In the above expression, fixing $z_1=z_2=z\in D_{\delta_1}$, we obtain
$$\begin{aligned}&\partial_{\nu_a}(A+ze^{i\alpha\pi})^{-1}(A+ze^{-i\alpha\pi})^{-1}\rho^{-1}f(x)\\
&=\sum_{n=0}^\infty \left[\sum_{k+j=n}(-e^{i\alpha\pi})^{k}(-e^{-i\alpha\pi})^{j}\right]\partial_{\nu_a}A^{-n-2}\rho^{-1}f(x)z^n\\
&=\sum_{n=0}^\infty \left[\sum_{k+j=n}e^{i\alpha\pi(k-j)}\right](-1)^n\partial_{\nu_a}A^{-n-2}\rho^{-1}f(x)z^n\\
&=\sum_{n=0}^\infty \left[\sum_{k=0}^ne^{i\alpha\pi (2k-n)}\right](-1)^n\partial_{\nu_a}A^{-n-2}\rho^{-1}f(x)z^n\\
&=\sum_{n=0}^\infty \left[\frac{1-e^{i2\alpha \pi(n+1)}}{1-e^{2i\alpha \pi}}\right]e^{-i\alpha\pi n}(-1)^n\partial_{\nu_a}A^{-n-2}\rho^{-1}f(x)z^n.\end{aligned}$$
Recalling that   $x\in \partial\Omega\setminus\mathcal G(f)$ if and only if $x\in\partial\Omega$ satisfies the following condition
$$\partial_{\nu_a}A^{-n-2}\rho^{-1}f(x)=0,\quad n\in\mathbb N\cup\{0\},\ \alpha(n+1)\not\in \mathbb N,$$
we deduce that
$$\left[\frac{1-e^{i2\alpha \pi(n+1)}}{1-e^{2i\alpha \pi}}\right]e^{-i\alpha\pi n}(-1)^n\partial_{\nu_a}A^{-n-2}\rho^{-1}f(x)=0,\quad
x\in \partial\Omega\setminus\mathcal G(f),\ n\in\mathbb N\cup\{0\}.$$
Therefore, we have
$$\partial_{\nu_a}(A+ze^{i\alpha\pi})^{-1}(A+ze^{-i\alpha\pi})^{-1}\rho^{-1}f(x)=0,\quad z\in D_{\delta_1},\ x\in \partial\Omega\setminus\mathcal G(f)$$
which implies that for all $x\in \partial\Omega\setminus\mathcal G(f)$ condition \eqref{cond3} is fulfilled. This implies that $\partial\Omega\setminus\mathcal G(f)\subset\partial\Omega\setminus\mathcal J(f)$ and we deduce $\mathcal J(f)\subset \mathcal G(f)$. In the same way,  from the above argumentation, we deduce that for all $x\in\partial\Omega$ the map $$z\mapsto \partial_{\nu_a}(A+ze^{i\alpha\pi})^{-1}(A+ze^{-i\alpha\pi})^{-1}\rho^{-1}f(x)$$ is holomorphic with respect to $z\in D_{\delta_1}$ and we have
$$\partial_{\nu_a}(A+ze^{i\alpha\pi})^{-1}(A+ze^{-i\alpha\pi})^{-1}\rho^{-1}f(x)=\sum_{n=0}^\infty \left[\frac{1-e^{i2\alpha \pi(n+1)}}{1-e^{2i\alpha \pi}}\right]e^{-i\alpha\pi n}(-1)^n\partial_{\nu_a}A^{-n-2}\rho^{-1}f(x)z^n.$$
Thus, fixing $x\in \partial\Omega\setminus\mathcal J(f)$, condition \eqref{cond3} implies
$$\partial_{\nu_a}(A+ze^{i\alpha\pi})^{-1}(A+ze^{-i\alpha\pi})^{-1}\rho^{-1}f(x)=0,\quad z\in D_{\delta_1}.$$
Then, it follows that
$$\begin{aligned}&n!\left[\frac{1-e^{i2\alpha \pi(n+1)}}{1-e^{2i\alpha \pi}}\right]e^{-i\alpha\pi n}(-1)^n\partial_{\nu_a}A^{-n-2}\rho^{-1}f(x)\\
&=\partial_{z}^n\left.\left[\partial_{\nu_a}(A+ze^{i\alpha\pi})^{-1}(A+ze^{-i\alpha\pi})^{-1}\rho^{-1}f(x)\right]\right|_{z=0}=0,\quad n\in\mathbb N\cup\{0\}.\end{aligned}$$
Therefore, we have
$$\left[\frac{1-e^{i2\alpha \pi(n+1)}}{1-e^{2i\alpha \pi}}\right]\partial_{\nu_a}A^{-n-2}\rho^{-1}f(x)=0, \quad x\in \partial\Omega\setminus\mathcal J(f),\ n\in\mathbb N\cup\{0\}$$
and we deduce that
$$\partial_{\nu_a}A^{-n-2}\rho^{-1}f(x)=0, \quad x\in \partial\Omega\setminus\mathcal J(f),\ n\in\mathbb N\cup\{0\},\ (n+1)\alpha\not\in\mathbb N.$$
From this identity, one can easily check that for any $x\in \partial\Omega\setminus\mathcal J(f)$ we have $x\in \partial\Omega\setminus\mathcal G(f)$ and combining this with the fact that $\partial\Omega\setminus\mathcal G(f)\subset\partial\Omega\setminus\mathcal J(f)$ we deduce that $\partial\Omega\setminus\mathcal G(f)=\partial\Omega\setminus\mathcal J(f)$. It follows that $\mathcal G(f)=\mathcal J(f)$.
\ \\

\textbf{Step 2.} In this step we will prove that $\mathcal G(f)$ is dense in $\partial\Omega$. In view of Step 1, it would be enough to show that $\mathcal J(f)$ is dense in $\partial\Omega$. We will show this result by contradiction. Assuming that  $\mathcal J(f)$ is not dense in $\partial\Omega$, we deduce that there exists an open not empty set $\Gamma$ of $\partial\Omega$ such that $\Gamma\subset\partial\Omega\setminus\mathcal J(f)$. Then, for all $x\in\Gamma$ there exists $\epsilon_x>0$ such that
\begin{equation}\label{cond33}  \partial_{\nu_a}(A+re^{i\alpha\pi})^{-1}(A+re^{-i\alpha\pi})^{-1}\rho^{-1}f(x)=0,\quad r\in(0,\epsilon_x) .\end{equation}
On the other hand, using the fact that $\rho^{-1}f\in D(A^{s})$, we deduce that for all $x\in\partial\Omega$, we have
$$\partial_{\nu_a}(A+z)^{-1}(A+e^{-i2\alpha\pi}z)^{-1}\rho^{-1}f(x)=\sum_{n=1}^\infty \frac{\sum_{k=1}^{m_n}\left\langle \rho^{-1}f,\phi_{n,k}\right\rangle\partial_{\nu_a}\phi_{n,k}(x)}{(\lambda_n+z)(\lambda_n+e^{-i2\alpha\pi}z)},\quad z\in\mathcal O$$
and the map $z\mapsto \partial_{\nu_a}(A+z)^{-1}(A+e^{-i2\alpha\pi}z)^{-1}\rho^{-1}f(x)$ is holomorphic with respect to $z\in\mathcal O$.
Thus, the condition \eqref{cond33} implies that
$$\partial_{\nu_a}(A+z)^{-1}(A+e^{-i2\alpha\pi}z)^{-1}\rho^{-1}f(x)\equiv0,\quad z\in \mathcal O,\ x\in\Gamma$$
and we deduce that \eqref{t1d} holds true. Therefore, repeating the arguments used at the end of the proof of Theorem \ref{t1}, we deduce that \eqref{t1d} implies that $f\equiv0$ which contradicts the fact that $f\not\equiv0$. This proves the density of the set $\mathcal G(f)=\mathcal J(f)$  in $\partial\Omega$.
\ \\
\textbf{Step 3.} In this step we will show that for any $x_0\in\mathcal G(f)$ the implication \eqref{t2b} holds true. For this purpose, we fix $x_0\in\mathcal G(f)$ and from   Step 1 we deduce that $x_0\in\mathcal G(f)=\mathcal J(f)$. Therefore, for all $\epsilon'>0$ there exists $r_0\in (0,\epsilon')$ such that \eqref{cond2} is fulfilled.
Let $u$ be the weak solution of problem \eqref{eq1}  satisfying
$$\partial_{\nu_a}u(x_0,t)=0,\quad t\in(T-\epsilon,T).$$
Fixing $v$ the weak solution of \eqref{eq2} and recalling that $v=u$ on $\Omega\times(0,T)$, we deduce that
\begin{equation}\label{t2h}\partial_{\nu_a}v(x_0,t)=0,\quad t\in(T-\epsilon,T).\end{equation}
Without loss of generality we may assume that $\delta=\epsilon$ with $\delta$ appearing in condition \eqref{source}.
 Then, applying Lemma \ref{l2} and the unique continuation of analytic functions we deduce that
$$\partial_{\nu_a}v(x_0,t)=0,\quad t\in(T-\epsilon,+\infty)$$
and fixing $h:=\R_+\ni t\mapsto \partial_{\nu_a}v(x_0,t)$ we deduce that supp$(h)\subset[0,T-\epsilon]$. Therefore, the Laplace transform $\hat{h}$ of $h$ admits an holomorphic extension to $\mathbb C$. Moreover, \eqref{l2a} implies that
$$h(t)=\int_0^t\sigma(s) \partial_{\nu_a} [S(t-s)\rho^{-1}f](x_0)ds,\quad t\in\R_+.$$
Applying the Laplace transform on both side of this identity and using \eqref{l2b}, we get
\begin{equation}\label{t2i}\hat{h}(p)=\hat{\sigma}(p)\sum_{n=1}^\infty \frac{\sum_{k=1}^{m_n}\left\langle \rho^{-1}f,\phi_{n,k}\right\rangle\partial_{\nu_a}\phi_{n,k}(x_0)}{\lambda_n+p^\alpha},\quad p\in\mathbb C_+.\end{equation}
We set $R>0$ and we consider the identity \eqref{t1b} for $p=Re^{i\theta}$ with $\theta\in(-\pi,\pi)$. In a similar way to the proof of Theorem \ref{t1}, sending $\theta\to\pm\pi$ and using the fact that $\hat{h}(p)$ and $\hat{\sigma}(p)$ are holomorphic with respect  to $p\in\mathbb C$, we obtain
\bel{t2aa}\hat{\sigma}(-R)\partial_{\nu_a}(A+R^\alpha e^{i\alpha\pi})^{-1}(A+R^\alpha e^{-i\alpha\pi})^{-1}\rho^{-1}f(x_0)= 0,\quad R\in \R_+.\ee
Combining \eqref{t2aa} with condition \eqref{cond2}, we deduce that for all $n\in\mathbb N$, there exists $r_n\in (0,2^{-n})$ such that
$$\partial_{\nu_a}(A+r_n^\alpha e^{i\alpha\pi})^{-1}(A+r_n^\alpha e^{-i\alpha\pi})^{-1}\rho^{-1}f(x_0)\neq0.$$
Then \eqref{t2aa} implies
\bel{t2ab}\hat{\sigma}(-r_n)=0,\quad n\in\mathbb N.\ee
It is clear that the sequence $(r_n)_{n\in\mathbb N}$ admits an accumulation at $0$ and using the fact that $\hat{\sigma}$ is holomorphic on $\mathbb C$, we deduce that $\hat{\sigma}\equiv 0$. Therefore, we have $\sigma\equiv0$. This proves that the implication \eqref{t1a} holds true.
\ \\

\textbf{Step 4.} In this step, we prove the last statement of Theorem \ref{t1}. For this purpose, we fix $k_1\in\mathbb N\cup\{0\}$ and $g\in D(A^{k_1+s})$ of constant sign. Then, we fix $f=\rho A^{k_1}g$. By eventually replacing $f$ by $-f$, we may assume without loss of generality that  $g\leq0$. We fix $k_2\in \mathbb N$ such that $k_2\geq k_1$ and $\alpha(k_2+1)\not\in\mathbb N$. According to the above discussion the proof will be completed if we prove that for any $x_0\in\partial\Omega$ we have $\partial_{\nu_a}A^{-2-k_2}\rho^{-1}f(x_0)\neq0$.

Set $x_0\in\partial\Omega$ and $w=A^{-2-k_2}\rho^{-1}f=A^{-2-k_2+k_1}g$. Using the fact that $\mathcal A A^{-1}g= \rho g\leq 0$ and $A^{-1}g|_{\partial\Omega}=0$, the maximum principle implies that $A^{-1}g\leq0$. In the same way by iteration, we can prove that $A^{-1-k_2+k_1}g\leq0$. Therefore, since $\mathcal Aw= \rho A^{-1-k_2+k_1}g\leq 0$ and $w|_{\partial\Omega}=0$,  the strong maximum principle (see e.g. \cite[Corollary 3.5]{GT}) implies that $w(x)<0=A^{-2-k_2+k_1}g(x_0),\   x\in\Omega$.
Thus,  the Hopf lemma (see \cite[Lemma 3.4]{GT}) implies that \begin{equation}\label{t2j}\partial_{\nu} A^{-2-k_2}\rho^{-1}f(x_0)=\partial_{\nu}w(x_0)>0.\end{equation}
On the other hand, recalling   that $w|_{\partial\Omega}\equiv0$ we deduce that $\nabla w(x_0)=(\partial_\nu w(x_0))\nu(x_0)$. Then, fixing the matrix $C=(a_{ij}(x_0))_{1\leq i,j\leq d}$, we get
$$\partial_{\nu_a} w(x_0)=[\nabla w(x_0)]\cdot [C\nu(x_0)]\\
=[(\partial_\nu w(x_0))\nu(x_0)]\cdot [C\nu(x_0)]\\
=(\partial_\nu w(x_0))\nu(x_0)^TC\nu(x_0).$$
Moreover, according to \eqref{ell} we have
$$\nu(x_0)^TC\nu(x_0)= \sum_{i,j=1}^d a_{i,j}(x_0) \nu_i(x_0) \nu_j(x_0) \geq c |\nu(x_0)|^2=c>0.$$
Combining this with \eqref{t2j}, we deduce that $\partial_{\nu_a} A^{-2-k_2}\rho^{-1}f(x_0)=\partial_{\nu_a} w(x_0)>0$. Therefore, we have $x_0\in\mathcal G(f)$ and it follows that $\mathcal G(f)=\partial\Omega$. Moreover, according to the Step 3 of the proof of this theorem, the implication \eqref{t2b} holds true for any $x_0\in\partial\Omega$.
This completes the proof of the last statement of the theorem and by the same way the proof of the theorem.

\section*{Aknowledgement}
The work of the first author is supported by the Grant PRG832 of the Estonian Research Council. The work of the  second  author is partially supported by the Agence Nationale de la Recherche (project MultiOnde) under  grant ANR-17-CE40-0029.

%
\end{document}